\documentclass[english,11pt]{article}
\usepackage[latin1]{inputenc}
\usepackage{amsmath,amsthm,amssymb}
\usepackage{graphics,epsfig,graphicx}

\def\proof{{\it Proof.}\ }
\def\endproof{\nolinebreak\hfill $\square$ \par\vskip3mm}

\def\eq#1{(\ref{#1})}

\def\neweq#1{\begin{equation}\label{#1}}
\def\endeq{\end{equation}}
\def\weak{\rightharpoonup}
\def\ep{\varepsilon}
\def\la{\lambda}

\def\RR{{\mathbb R} }

\def\di{\displaystyle}
\def\ri{\rightarrow}

\def\ii{\^\i }
\newtheorem{theorem}{Theorem}[section]
\newtheorem{lem}{Lemma}[section]

\title{\sc Multiparameter bifurcation and asymptotics for the singular
Lane-Emden-Fowler equation with convection term}
\author{Marius GHERGU and
Vicen\c tiu R\u ADULESCU\thanks{Correspondence address: Vicen\c
tiu R\u adulescu, Department of Mathematics, University of
Craiova, 200585 Craiova, Romania, fax: +40-251.41.16.88. E-mail:
{\tt radulescu@inf.ucv.ro}}\\
\small Department of Mathematics, University of Craiova, 200585
Craiova, Romania}
\date{}

\begin{document}
\baselineskip16pt \maketitle
\renewcommand{\theequation}{\arabic{section}.\arabic{equation}}
\catcode`@=11 \@addtoreset{equation}{section} \catcode`@=12
\begin{abstract}
We establish some bifurcation results for the boundary value
problem $-\Delta u=g(u)+\la|\nabla u|^p+\mu f(x,u)$ in $\Omega,$
$u>0$ in $\Omega,$ $u=0$ on $\partial\Omega,$ where $\Omega$ is a
smooth bounded domain in $\RR^N,$ $\la,\mu\geq 0,$ $0<p\leq 2,$
$f$ is nondecreasing with respect to the second variable, and $g$
is unbounded around the origin. The asymptotic behaviour of the
solution around the bifurcation point is also established,
provided $g(u)$ behaves like $u^{-\alpha}$ around the origin, for
some $0<\alpha<1$. Our approach relies on finding explicit sub-
and super-solutions combined with various techniques related to
the maximum principle for elliptic equations. The analysis we
develop in this paper shows the key role played by the convection
term $|\nabla u|^p$.

\noindent{\bf Key words}: singular elliptic equation, sub-linear
perturbation,
 bifurcation problem,
maximum principle.

\noindent{\bf 2000 Mathematics Subject Classification}: 35A20,
35B25, 35B50, 35J60, 58J55, 58K55.
\end{abstract}

\section{Introduction and the main results}
In his recent monograph \cite{kiel}, H.~Kielh\"ofer synthetizes
the role of bifurcation problems in Applied Mathematics: {\it
Bifurcation Theory attempts to explain various phenomena that have
been discovered and described in Natural Sciences over the
centuries. The buckling of the Euler rod, the appearance of Taylor
vortices, and the onset of oscillations in an electric circuit,
for instance, all have a common cause: a specific physical
parameter crosses a threshold, and that event forces the system to
the organization of a new state that differs considerably from
that observed before}.

In the present paper we continue the bifurcation analysis
developed in our previous works \cite{gr, grcras} (see also
\cite{crg}) for a large class of semilinear elliptic equations
with singular nonlinearity and Dirichlet boundary condition. Such
problems arise in the study of non-Newtonian fluids, boundary
layer phenomena for viscous fluids, chemical heterogeneous
catalysts, as well as in the theory of heat conduction in
electrically conducting materials.
 The main feature of this paper
is the presence of the convection term $|\nabla u|^p$.

Let $\Omega\subset \RR^N$ $(N\geq 2)$ be a bounded domain with a
smooth boundary. We are concerned in this paper with singular
elliptic problems of the following type
\neweq{P}
\left\{\begin{tabular}{ll}
$-\Delta u=g(u)+\la|\nabla u|^p+\mu f(x,u)$ \quad & $\mbox{\rm in}\ \Omega,$\\
$u>0$ \quad & $\mbox{\rm in}\ \Omega,$\\
$u=0$ \quad & $\mbox{\rm on}\ \partial\Omega,$\\
\end{tabular} \right.
\endeq
where $0<p\leq 2$ and $\la,\mu\geq 0.$ As remarked in \cite{leray,
kazdan}, the requirement that the nonlinearity grows at most
quadratically in $|\nabla u|$ is natural in order to apply the
maximum principle.

Throughout this paper we suppose that
$f:\overline{\Omega}\times[0,\infty)\rightarrow[0,\infty)$ is a
H\"{o}lder continuous function which is nondecreasing with respect
to the second variable and is positive on $\overline\Omega\times
(0,\infty).$ We assume that $g:(0,\infty)\ri(0,\infty)$ is a
H\"{o}lder continuous function which is nonincreasing and
$\lim_{s\searrow 0}g(s)=+\infty.$

Problems of this type arise in the study of guided modes of an
electromagnetic field in a nonlinear medium, satisfying adequate
constitutive hypotheses. The following two examples illustrate
situations of this type: (i) if $f(u)=u^3(1+\gamma u^2)^{-1}$
($\gamma>0$) then problem \eq{P} describes the variation of the
dielectric constant of gas vapors where a laser beam propagates
(see \cite{stuart1, stuart2}); (ii) nonlinearities of the type
$f(u)=(1-e^{-\gamma u^2})u$ arise in the context of laser beams in
plasmas (see~\cite{sulem}). If $f(u)=e^{u/(1+\varepsilon u)}$
($\ep>0$) then the corresponding equation describes the
temperature dependence of the reaction rate for exothermic
reactions obeying the simple Arrhenius rate law in circumstances
in which the heat flow is purely conductive (see \cite{chim1,
chim2}). In this context the parameter $\ep$ is a dimensionless
ambient temperature and the parameter $\lambda$ is a dimensionless
heat evolution rate. The corresponding equation
$$-\Delta u=g(u)+\la|\nabla u|^p+\mu e^{u/(1+\varepsilon
u)}\qquad\mbox{in}\ \,\Omega$$ represents heat balance with
reactant consumption ignored, where $u$ is a dimensionless
temperature excess. The Dirichlet boundary condition $u=0$ on
$\partial\Omega$ is an isothermal condition and, in this case, it
describes the exchange of heat at the surface of the reactant by
Newtonian cooling.

Our general setting includes some simple prototype models from
boundary-layer theory of viscous fluids (see \cite{wong}). If
$\la=0$ and $\mu=0,$ \eq{P} is called the Lane-Emden-Fowler
equation. Problems of this type, as well as the associated
evolution equations, describe naturally  certain physical
phenomena. For example, super-diffusivity equations of this type
have been proposed by de~Gennes \cite{gennes} as a model for long
range Van der Waals interactions in thin films spreading on solid
surfaces. This equation also appears in the study of cellular
automata and interacting particle systems with self-organized
criticality (see \cite{chayes}), as well as to describe the flow
over an impermeable plate (see \cite{cn1, cn}).

 Our aim in this paper is to describe the influence of the
gradient term in problem \eq{P}.

 Many papers have been devoted to the case $\la=0,$ where the
problem \eq{P} becomes
\neweq{la}
\left\{\begin{tabular}{ll}
$-\Delta u=g(u)+\mu f(x,u)$ \quad & $\mbox{\rm in}\ \Omega,$\\
$u>0$ \quad & $\mbox{\rm in}\  \Omega,$\\
$u=0$ \quad & $\mbox{\rm on}\ \partial\Omega,$\\
\end{tabular} \right.
\endeq
If $\mu=0,$ then \eq{la} has a unique solution (see \cite{crt,
lm1}). When $\mu>0,$ the study of \eq{la} emphasizes the role
played by the nonlinear term $f(x,u).$  For instance, if one of
the following assumptions are fulfilled

$(f1)\;$ there exists $c>0$ such that $f(x,s)\geq cs$ for all
$(x,s)\in \overline\Omega\times [0,\infty);$

$(f2)\;$ the mapping $\di (0,\infty)\ni s\longmapsto
\frac{f(x,s)}{s}$
is nondecreasing for all $x\in\overline\Omega,$\\
then problem \eq{la} has solutions only if $\mu>0$ is small enough
(see \cite{cp}). In turn, when $f$ satisfies the following
assumptions

$\di(f3)\;$ the mapping $\di (0,\infty)\ni
s\longmapsto\frac{f(x,s)}{s}\quad\mbox{is nonincreasing for
all}\;\, x\in\overline{\Omega};$

$\di(f4)\; \lim_{s\rightarrow\infty}\frac{f(x,s)}{s}=0,\;\;
\mbox{uniformly for}\;\,x\in\overline{\Omega},$\\
then problem \eq{la} has at least one solutions for all $\mu>0$
(see \cite{crg, cp, gr, sy1} and the references therein). The same
assumptions will be used in the study of \eq{P}.

If $\la>0,$ the following problem was treated in Zhang and Yu
\cite{zy}
\neweq{aunu}
\left\{\begin{tabular}{ll} $\di -\Delta
u=\frac{1}{u^\alpha}+\la|\nabla u|^p+\sigma$ \quad &
$\mbox{\rm in}\ \Omega,$\\
$u>0$ \quad & $\mbox{\rm in}\ \Omega,$\\
$u=0$ \quad & $\mbox{\rm on}\ \partial\Omega,$\\
\end{tabular} \right.
\endeq
where $\la$, $\sigma \geq 0,$ $\alpha>0$, and $p\in (0,2]$.
 By using the change of variable $v=e^{\la u}-1$ in the
case $p=2,$ it is proved in \cite{zy} that problem \eq{aunu} has
classical solutions if $\la\sigma<\la_1,$ where $\la_1$ is the
first eigenvalue of $-\Delta$ in $H^1_0(\Omega).$ This will be
used to deduce the existence and nonexistence in the case $0<p<2.$

 If $f(x,u)$
depends on $u,$ the above change of variable does not preserve the
sublinearity condition $(f3)-(f4)$ and the monotony of the
nonlinear term $g$ in problem \eq{P}. In turn, if $f(x,u)$ does
not depend on $u$ and $p=2,$ this method successfully applies to
our study and we will be able to give a complete characterization
of \eq{P} (see Theorem \ref{th4} below).

Due to the singular term $g(u)$ in $(P_\la),$ we cannot expect to
have solutions in $C^2(\overline\Omega).$ As it was pointed out in
\cite{zy}, if $\alpha>1$ then the solution of \eq{aunu} is {\bf
not} in $C^1(\overline\Omega).$ We are seeking in this paper
classical solutions of $(P_\la),$ that is, solutions $u\in
C^2(\Omega)\cap C(\overline\Omega)$ that verify \eq{P}.
\medskip

By the monotony of $g,$ there exists
$$\di a=\lim_{s\ri\infty}g(s)\in[0,\infty).$$

The first result concerns the case $\la=1$ and $1<p\leq 2.$ In the
statement of the following result we do not need assumptions
$(f1)-(f4);$ we just require that $f$ is a H\"{o}lder continuous
function which is nondecreasing with respect to the second
variable and is positive on $\overline \Omega\times(0,\infty).$
\begin{theorem}\label{th1} Assume $\la=1$ and $1<p\leq 2.$\\
{\rm (i) } If $p=2$ and $a\geq \la_1,$ then \eq{P} has no solutions;\\
{\rm (ii) } If $p=2$ and $a<\la_1$ or $1<p<2,$ then there exists
$\mu^*>0$ such that \eq{P} has at least one classical solution for
$\mu<\mu^*$ and no solutions exist if $\;\mu>\mu^*.$
\end{theorem}

If $\la=1$ and $0<p\leq 1$ the study of existence is close related
to the asymptotic behaviour of the nonlinear term $f(x,u).$ In
this case we prove

\begin{theorem}\label{th2} Assume $\la=1$ and $0<p\leq 1.$\\
{\rm (i) } If $f$ satisfies $(f1)$ or $(f2),$ then there exists
$\mu^*>0$ such that \eq{P} has at least one classical solution for
$\mu<\mu^*$ and no solutions exist if $\;\mu>\mu^*;$\\
{\rm (ii) } If $\;0<p<1$ and $f$ satisfies $(f3)-(f4),$ then
\eq{P} has at least one solution for all $\mu\geq 0.$
\end{theorem}

 Next we are concerned with the case $\mu=1.$ Our result is the
following

\begin{theorem}\label{th3} Assume $\mu=1$ and $f$ satisfies
assumptions $(f3)$ and $(f4).$ Then the following properties hold true.\\
{\rm (i) } If $\;0<p<1,$ then \eq{P} has at least one classical
solution for all $\la\geq 0$;\\
{\rm (ii) } If $\;1\leq p\leq 2,$ then there exists
$\la^*\in(0,\infty]$ such that \eq{P} has at least one classical
solution for $\la<\la^*$ and no solution exists if $\;\la>\la^*.$
Moreover, if $\,1<p\leq 2,$ then $\la^*$ is finite.
\end{theorem}
Related to the above result we raise the following {\bf open
problem:} if $p=1$ and $\mu=1,$ is $\la^*$ a finite number?

Theorem \ref{th3} shows the importance of the convection term
$\lambda |\nabla u|^p$ in \eq{P}. Indeed, according to
\cite[Theorem~1.3]{gr} and for any $\mu>0$, the boundary value
problem
\neweq{PQ}
\left\{\begin{tabular}{ll}
$-\Delta u=u^{-\alpha}+\la|\nabla u|^p+\mu u^\beta$ \quad & $\mbox{\rm in}\ \Omega,$\\
$u>0$ \quad & $\mbox{\rm in}\ \Omega,$\\
$u=0$ \quad & $\mbox{\rm on}\ \partial\Omega$\\
\end{tabular} \right.
\endeq
has a unique solution, provided $\lambda=0$, $\alpha$, $\beta\in
(0,1)$. The above theorem shows that if $\la$ is not necessarily
0, then the following situations may occur : (i) problem \eq{PQ}
has solutions if $p\in (0,1)$ and for all $\lambda\geq 0$; (ii) if
$p\in (1,2)$ then there exists $\la^*>0$ such that problem \eq{PQ}
has a solution for any $\la<\la^*$ and no solution exists if
$\la>\la^*.$

To see the dependence between $\la$ and $\mu$ in \eq{P}, we
consider the special case $f\equiv 1$ and $p=2.$ In this case we
can say more about the problem \eq{P}. More precisely we have

\begin{theorem}\label{th4}
Assume that $p=2$ and $f\equiv 1.$\\
{\rm (i) } The problem \eq{P} has solution if and only if $\la(a+\mu)<\la_1;$\\
{\rm (ii) } Assume $\mu>0$ is fixed, $g$ is decreasing and let
$\di \la^*=\frac{\la_1}{a+\mu}.$ Then \eq{P} has a unique solution
$u_\la$ for all $\la<\la^*$ and the sequence $(u_\la)_{\la<\la^*}$
is increasing with respect to $\la.$\\
Moreover, if $\,\di \limsup_{s\searrow 0}s^\alpha g(s)<+\infty,$
for some $\alpha\in(0,1),$ then the sequence of solutions
$(u_\la)_{0<\la<\la^*}$ has the following properties

\qquad{\rm (ii1)} For all $0<\la<\la^*$ there exist two positive
constants $c_1,c_2$ depending on $\la$ such that $c_1\,{\rm
dist}(x,\partial \Omega)\leq u_\la\leq c_2\,{\rm dist}(x,\partial
\Omega)$ in $\Omega;$

$\qquad${\rm (ii2)} $u_\la\in C^{1,1-\alpha}(\overline\Omega)\cap
C^2(\Omega);$

$\qquad${\rm (ii3)} $u_{\la}\longrightarrow +\infty$ as
$\la\nearrow \la^*$, uniformly on compact subsets of $\Omega.$
\end{theorem}

The assumption $\di \limsup_{s\searrow 0}s^\alpha g(s)<+\infty,$
for some $\alpha\in(0,1)$, has been used in \cite{gr} and it
implies the following Keller-Osserman-type growth condition around
the origin
\neweq{KellerOser}\int_0^1\left(\int_0^tg(s)ds\right)^{-1/2}dt<+\infty.\endeq
As proved by B\'enilan, Brezis and Crandall in \cite{brebc},
condition \eq{KellerOser} is equivalent to the {\it property of
compact support}, that is, for any $h\in L^1(\RR^N)$ with compact
support, there exists a unique $u\in W^{1,1}(\RR^N)$ with compact
support such that $\Delta u\in L^1(\RR^N)$ and
$$-\Delta u=g(u)+h\qquad\mbox{a.e. in}\ \RR^N.$$

The situations described in Theorem~\ref{th4} are depicted in the
following bifurcation diagrams. Case~1 (resp., Case~2) corresponds
to (i) and $a=0$ (resp., $a>0$), while Case~3 is related to (ii),
$\lambda>0$ and $\mu=\mbox{fixed}$.
\begin{figure}[ht]
\begin{center}
\epsfig{file=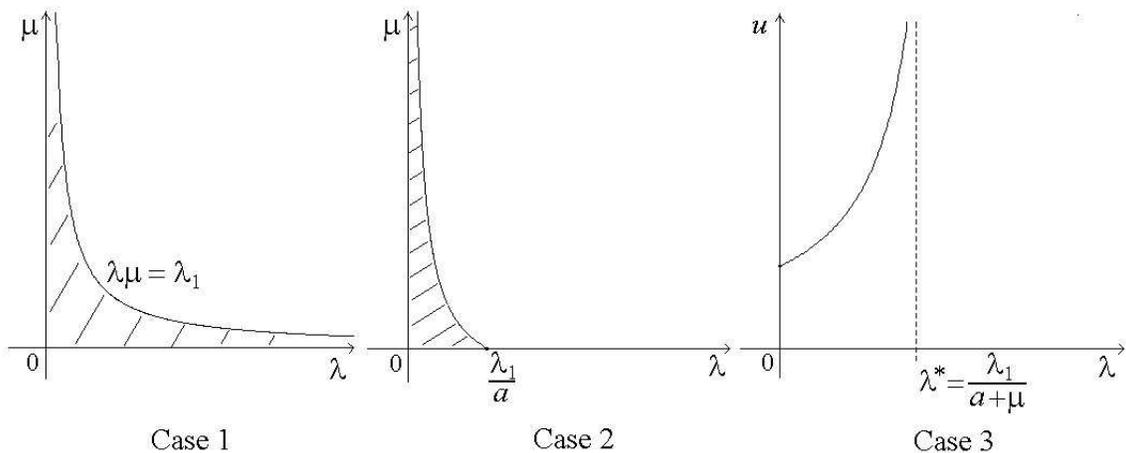,height=6.0cm} \caption {Bifurcation
diagrams}
\end{center}
%\vskip -0.5cm
\end{figure}

As regards the uniqueness of the solutions to problem \eq{P}, we
may say that this does not seem to be a feature easy to achieve.
Only when $f(x,u)$ is constant in $u$ we can use classical methods
in order to prove the uniqueness. It is worth pointing out here
that the uniqueness of the solution is a delicate issue even for
the simpler problem \eq{la}. We have showed in \cite{gr} that when
$f$ fulfills $(f3)-(f4)$ and $g$ satisfies the same growth
condition as in Theorem \ref{th4} then, if \eq{la} has a solution,
it will follows that this solution is unique. On the other hand,
if $f$ satisfies $(f2),$ the uniqueness generally does not occur.
In that sense we refer the interested reader to Haitao
\cite{haitao}. In the case $f(x,u)=u^q,$ $g(u)=u^{-\gamma},$
$0<\gamma<\frac{1}{N}$ and $1<q<\frac{N+2}{N-2},$ we learn from
\cite{haitao} that problem \eq{la} has at least two classical
solutions provided $\mu$ belongs to a certain range.

\smallskip
Our approach relies on finding of appropriate sub- and
super-solutions of \eq{P}. This will allows us to enlarge the
study of bifurcation to a class of problems more generally to that
studied in \cite{zy}. However, neither the method used in Zhang
and Yu \cite{zy}, nor our method gives a precise answer if $\la^*$
is finite or not in the case $p=1$ and $\mu=1.$

In the next Section we state some auxiliary results which will be
used in the proofs of the above Theorems. This will be done in
Sections 3, 4, 5 and 6.

\section{Auxiliary results}
Let $\varphi_1$ be the normalized positive eigenfunction
corresponding to the first eigenvalue $\la_1$ of the problem
$$ \left\{\begin{tabular}{ll}
$-\Delta u=\la u$ \quad & $\mbox{\rm in}\ \Omega,$\\
$u=0$ \quad & $\mbox{\rm on}\ \partial\Omega,$\\
\end{tabular} \right.$$
As it is well known $\la_1>0,$ $\varphi_1\in
C^2({\overline{\Omega}})$ and
\neweq{fi}
C_1\,{\rm dist}(x,\partial \Omega)\leq\varphi_1\leq C_2\,{\rm
dist}(x,\partial \Omega)\quad\mbox{ in }\;\Omega,
\endeq
for some positive constants $C_1,C_2>0.$ From the characterization
of $\la_1$ and $\varphi_1$ we state the following elementary
result. For the convenience of the reader we shall give a complete
proof.
\begin{lem}\label{l01} Let $F:\overline{\Omega}\times(0,\infty)\rightarrow\RR$ be a
continuous function such that $F(x,s)\geq \la_1s+b$ for some $b>0$
and for all $(x,s)\in \overline\Omega\times (0,\infty).$ Then the
problem
\neweq{la1}
\left\{\begin{tabular}{ll}
$-\Delta u=F(x,u)$ \quad & $\mbox{\rm in}\ \Omega,$\\
$u> 0$ \quad & $\mbox{\rm in}\ \Omega,$\\
$u=0$\quad & $\mbox{\rm on}\ \partial\Omega,$\\
\end{tabular} \right.\endeq
has no solutions.
\end{lem}
\noindent \proof By contradiction, suppose that \eq{la1} admits a
solution. This will provide a super-solution of the problem
\neweq{la2}
\left\{\begin{tabular}{ll}
$-\Delta u=\la_1 u+b$ \quad & $\mbox{\rm in}\ \Omega,$\\
$u> 0$ \quad & $\mbox{\rm in}\ \Omega,$\\
$u=0$\quad & $\mbox{\rm on}\ \partial\Omega,$\\
\end{tabular} \right.\endeq
Since 0 is a sub-solution, by the sub and super-solution method
and classical regularity theory it follows that \eq{la1} has a
solution $u\in C^2(\overline\Omega).$ Multiplying by $\varphi_1$
in \eq{la2} and then integrating over $\Omega,$ we get
$$\di -\int_{\Omega}\varphi_1\Delta u =\la_1\int_\Omega \varphi_1u+b\int_{\Omega}\varphi_1,$$
that is $\di \la_1\int_{\Omega}\varphi_1 u =\la_1\int_\Omega
\varphi_1u+b\int_{\Omega}\varphi_1,$ which implies $\di
\int_\Omega\varphi_1=0.$ This is clearly a contradiction since
$\varphi_1>0$ in $\Omega.$ Hence \eq{la1} has no solutions.
\endproof
The growth of $\varphi_1$ is prescribed in the following result.
\begin{lem}\label{l1}\emph{(see \cite{lm1}).}
$\di \int_{\Omega}\varphi_1^{-s}dx<+\infty$ if and only if $s<1.$
\end{lem}
Basic in the study of the existence is the following Lemma.
\begin{lem}\label{l2}\emph{(see \cite{sy1}).}
Let $F:\overline{\Omega}\times(0,\infty)\rightarrow\RR$ be a
H\"{o}lder continuous function on each compact subset of
$\overline{\Omega}\times(0,\infty)$ which satisfies
\medskip

\noindent $\di(F1)\qquad \limsup _{s\rightarrow +\infty}
\left(s^{-1}\max_{x\in\overline{\Omega}}F(x,s)\right)<\la_1;$
\smallskip

\noindent $(F2)\qquad$ for each $t>0,$ there exists a constant
$D(t)>0$ such that
$$F(x,r)-F(x,s)\geq -D(t)(r-s),\quad\mbox{for}\;\;x\in {\overline{\Omega}}\,
\;\;\mbox{and}\;\;\;r\geq s\geq t;$$
\smallskip

\noindent $(F3)\qquad$ there exist a $\eta_0>0$ and an open subset
$\Omega_0\subset\Omega$ such that
$$\di\min_{x\in\overline{\Omega}}F(x,s)\geq 0\quad\mbox{for}\;s\in(0,\eta_0),$$
and
$$\di \lim_{s\searrow 0}\frac{F(x,s)}{s}=+\infty\quad\mbox{uniformly for}
\;\,x\in\Omega_0.$$

Then the problem
\begin{equation}\label{doidoi}
 \left\{\begin{tabular}{ll}
$-\Delta u=F(x,u)$ \quad & $\mbox{\rm in}\ \Omega,$\\
$u> 0$ \quad & $\mbox{\rm in}\ \Omega,$\\
$u=0$\quad & $\mbox{\rm on}\ \partial\Omega,$\\
\end{tabular} \right.
\end{equation}
has at least one solution.
\end{lem}

According to Lemma \ref{l2}, there exists $\zeta\in C^2({\overline
\Omega})$ a solution of the problem
\neweq{treidoi}
\left\{\begin{tabular}{ll}
$\di -\Delta \zeta=g(\zeta)$&$\mbox {\rm in}\ \Omega,$\\
$\di \zeta>0$&$\mbox {\rm in}\ \Omega,$\\
$\di \zeta=0$&$\mbox{\rm on}\ \partial\Omega.$\\
\end{tabular}\right.\endeq
Clearly $\zeta$ is a sub-solution of \eq{P} for all $\la\geq 0.$
It is worth pointing out here that the sub-super solution method
still works for the problem \eq{P}. With the same proof as in
Zhang and Yu \cite[Lemmma 2.8]{zy} that goes back to the
pioneering work of Amann \cite{a} we state the following result.
\begin{lem}\label{sub-sup}
Let $\la,\mu\geq 0.$ If \eq{P} has a super-solution $\overline
u\in C^2(\Omega) \cap C(\overline \Omega)$ such that
$\zeta\leq\overline u$ in $\Omega,$ then \eq{P} has at least a
solution.
\end{lem}
Another difficulty in the treatment of \eq{P} is the lack of the
usual maximal principle. The following result which is due to Shi
and Yao \cite{sy1} gives a comparison principle that applies to
singular elliptic equations.
\begin{lem}\label{l3}\emph{(see \cite{sy1}).}
Let $F:\overline{\Omega}\times[0,\infty)\rightarrow\RR$ be a
continuous function such that the mapping $(0,\infty)\ni
s\rightarrow\frac{F(x,s)}{s}$ is strictly decreasing at each
$x\in\Omega.$ Assume $v,w\in C^2(\Omega)\cap
C({\overline{\Omega}})$ and

$(a)\quad \Delta w+F(x,w)\leq 0\leq \Delta v+F(x,v)$ in $\Omega;$

$(b)\quad v,w>0$ in $\Omega$ and $v\leq w$ on $\partial\Omega;$

$(c)\quad\Delta v\in L^1(\Omega).$\\
Then $v\leq w$ in $\Omega.$
\end{lem}
\begin{lem}\label{l4}\emph{(see \cite{ap}).} If $\;p>1,$ then there exists
a real number $\bar \sigma>0$ such that the problem
\begin{equation}\label {doitrei}
 \left\{\begin{tabular}{ll}
$-\Delta u=|\nabla u|^p+\sigma$ \quad & $\mbox{\rm in}\ \Omega,$\\
$u=0$ \quad & $\mbox{\rm on}\ \partial\Omega,$\\
\end{tabular} \right.
\end{equation}
has no solutions for $\sigma>\bar\sigma.$
\end{lem}

\section{Proof of Theorem \ref{th1}}
We start with the following more general result.
\begin{lem}\label{lp2}
Let $F:\overline\Omega\times(0,\infty)\ri[0,\infty)$ and
$G:(0,\infty)\ri(0,\infty)$ be two H\"{o}lder continuous functions
that verify

$(A1)\;$ $F(x,s)>0,$ for all
$\;(x,s)\in\overline\Omega\times(0,\infty);$

$(A2)\;$ The mapping $[0,\infty)\ni s\longmapsto F(x,s)$ is
nondecreasing for all $x\in\overline\Omega;$

$(A3)\;$ $G$ is nonincreasing and $\;\lim_{s\searrow
0}G(s)=+\infty.$

Assume that $\tau>0$ is a positive real number. Then the following holds.\\
{\rm (i) } If $\;\tau \lim_{s\ri\infty}G(s)\geq \la_1,$ then the
problem
\neweq{trunu}
\left\{\begin{tabular}{ll}
$\di -\Delta u=G(u)+\tau|\nabla u|^2+\mu F(x,u)$&$\mbox{\rm in}\ \Omega,$\\
$\di u>0$&$\mbox{\rm in}\ \Omega,$\\
$\di u=0$&$\mbox{\rm on}\ \partial\Omega,$\\
\end{tabular}\right.
\endeq
has no solutions.\\
{\rm (ii) } If $\;\tau\lim_{s\ri\infty}G(s)<\la_1,$ then there
exists $\bar \mu>0$ such that the problem \eq{trunu} has at least
one solution for all $0\leq \mu<\bar\mu.$
\end{lem}
\proof (i) With the change of variable $\di v= e^{\tau u}-1,$ the
problem \eq{trunu} takes the form
\neweq{trdoi}
\left\{\begin{tabular}{ll}
$-\Delta v=\Psi_\mu(x,u)$\quad & $\mbox{\rm in}\ \Omega,$\\
$v>0$ \quad & $\mbox{\rm in}\ \Omega,$\\
$v=0$ \quad & $\mbox{\rm on}\ \partial\Omega,$\\
\end{tabular} \right. \endeq
where
$$\di \Psi_\mu(x,s)=\tau(s+1)G\left(\frac{1}{\tau}
\ln(s+1)\right)+\mu \tau(s+1)
F\left(x,\frac{1}{\tau}\ln(s+1)\right),$$
for all $(x,s)\in\overline\Omega\times (0,\infty).$\\
Taking into account the fact that $G$ is nonincreasing and
$\tau\lim_{s\ri\infty}G(s)\geq \la_1,$ we get
$$\Psi_\mu(x,s)\geq \la_1(s+1)\quad\mbox{ in }\;
\overline\Omega\times(0,\infty), \mbox{ for all }\; \mu\geq 0.$$
By Lemma \ref{l01} we conclude that \eq{trdoi} has no solutions.
Hence \eq{trunu} has no solutions.\\
(ii) Since
$$\di \lim _{s\rightarrow+\infty}
\frac{\tau(s+1)G\left(\frac{1}{\tau}\ln(s+1)\right)+1}{s}<\la_1$$
and
$$\di \lim _{s\searrow 0}
\frac{\tau(s+1)G\left(\frac{1}{\tau}\ln(s+1)\right)+1}{s}=+\infty,$$
we deduce that the mapping $(0,\infty)\ni s\longmapsto
\tau(s+1)G\left(\frac{1}{\tau}\ln(s+1)\right)+1$ fulfills the
hypotheses in Lemma \ref{l2}. According to this one, there exists
$\overline v\in C^2(\Omega)\cap C(\overline\Omega)$ a solution of
the problem
$$\left\{\begin{tabular}{ll}
$-\Delta v=\di
\tau(v+1)G\left(\frac{1}{\tau}\ln(v+1)\right)+1$\quad & $
\mbox{\rm in}\ \Omega,$\\
$v>0$ \quad & $\mbox{\rm in}\ \Omega,$\\
$v=0$ \quad & $\mbox{\rm in}\ \partial\Omega.$\\
\end{tabular} \right.$$
Define
$$\di \bar \mu:=\frac{1}{\tau (\|\overline v\|_{\infty}+1)}\cdot
\frac{1}{\di \max\limits_{x\in\overline\Omega}
F\left(x,\frac{1}{\tau}\ln(\|\overline v\|_{\infty}+1)\right)}.$$
It follows that $\overline v$ is a super-solution of \eq{trdoi}
for all
$0\leq \mu<\bar \mu.$\\
Next we provide a sub-solution $\underline v$ of \eq{trdoi} such
that $\underline v\leq \overline v$ in $\Omega.$ To this aim, we
apply Lemma \ref{l2} to get that there exists $\underline v\in
C^2(\Omega)\cap C(\overline\Omega)$ a solution of the problem
$$\left\{\begin{tabular}{ll}
$-\Delta v=\di \tau G\left(\frac{1}{\tau}\ln(v+1)\right)$\quad &
$\mbox{\rm in}\ \Omega,$\\
$v>0$ \quad & $\mbox{\rm in}\ \Omega,$\\
$v=0$ \quad & $\mbox{\rm on}\ \partial\Omega.$\\
\end{tabular} \right.$$
Clearly, $\underline v$ is a sub-solution of \eq{trdoi} for all
$0\leq \mu<\bar\mu.$ Let us prove now that $\underline v\leq
\overline v$ in $\Omega.$ Assuming the contrary, it follows that
$\max_{x\in\overline\Omega} \{\underline v-\overline v\}>0$ is
achieved in $\Omega.$ At that point, say $x_0,$ we have
$$\begin{aligned}
\di 0&\di \leq -\Delta (\underline v-\overline v)(x_0)\\
&\di \leq \tau \left[G\left(\frac{1}{\tau}\ln(\underline
v(x_0)+1)\right)-
G\left(\frac{1}{\tau}\ln(\overline v(x_0)+1)\right)\right]-1<0,\\
\end{aligned}$$
which is a contradiction. Thus, $\underline v\leq \overline v$ in
$\Omega.$ We have proved that $(\underline v,\overline v)$ is an
ordered pair of sub-super solutions of \eq{trdoi} provided $0\leq
\mu<\bar \mu.$ It follows that \eq{trunu} has at least one
classical solution for all $0\leq \mu<\bar \mu$ and the proof of
Lemma \ref{lp2} is now complete.
\endproof
\medskip

\noindent{\bf Proof of Theorem \ref{th1} completed}. According to
Lemma \ref{lp2}(i) we deduce that \eq{P} has no solutions if $p=2$
and $a\geq \la_1.$ Furthermore, if $p=2$ and $a<\la_1,$ in view of
Lemma \ref{lp2}(ii), we deduce that \eq{P} has at least one
classical solution if $\mu$ is small enough. Assume now $1<p<2$
and let us fix $C>0$ such that
\neweq{CC}
aC^{p/2}+C^{p-1}<\la_1.
\endeq
Define
$$\psi:[0,\infty)\ri[0,\infty),\quad \psi(s)=\frac{s^p}{s^2+C}.$$
A careful examination reveals the fact that $\psi$ attains its
maximum at $\bar s=\left(\frac{Cp}{2-p}\right)^{2-p}.$ Hence
$$\di \psi(s)\leq \psi(\bar s)=\frac{p^{p/2}(2-p)^{(2-p)/2}}{2C^{1-p/2}},
\quad\mbox{ for all }\;s\geq 0.$$ By the classical Young's
inequality we deduce
$$\di p^{p/2}(2-p)^{(2-p)/2}\leq 2,$$
which yields $\psi(s)\leq C^{p/2-1},$ for all $s\geq 0.$ Thus, we
have proved
\neweq{inq}
\di s^p\leq C^{p/2}s^2+C^{p/2-1},\quad \mbox{ for all }\;s\geq 0.
\endeq
Consider the problem
\neweq{trtrei}
\left\{\begin{tabular}{ll}
$-\Delta u=g(u)+C^{p/2-1}+C^{p/2}|\nabla u|^2+\mu f(x,u)$&$\mbox {in }\Omega$\\
$u>0$ \quad & $\mbox{\rm in}\ \Omega,$\\
$u=0$ \quad & $\mbox{\rm on}\ \partial\Omega,$\\
\end{tabular} \right. \endeq
By virtue of \eq{inq}, any solution of \eq{trtrei} is a
super-solution of \eq{P}.

 Using now \eq{CC} we get
$$\lim_{s\ri\infty}C^{p/2}(g(u)+C^{p/2-1})<\la_1.$$
The above relation enables us to apply Lemma \ref{lp2}(ii) with
$G(s)=g(s)+C^{p/2-1}$ and $\tau=C^{p/2}.$ It follows that there
exists $\bar\mu>0$ such that \eq{trtrei} has at least a solution
$u.$ With a similar argument to that used in the proof of Lemma
\ref{lp2}, we obtain $\zeta \leq u$ in $\Omega,$ where $\zeta$ is
defined in \eq{treidoi}. By Lemma \ref{sub-sup} we get that \eq{P}
has at least one solution if $0\leq \mu<\bar\mu.$
\medskip

We have proved that \eq{P} has at least one classical solution for
both cases $p=2$ and $a<\la_1$ or $1<p<2,$ provided $\mu$ is
nonnegative small enough. Define next
$$A=\{\mu\geq 0;\;\mbox{ problem }\eq{P}\mbox{ has at least one solution}\}.$$
The above arguments implies that $A$ is nonempty. Let $\mu^*=\sup
A.$ We first show that $[0,\mu^*)\subseteq A.$ For this purpose,
let $\mu_1\in A$ and $0\leq \mu_2<\mu_1.$ If $u_{\mu_1}$ is a
solution of \eq{P} with $\mu=\mu_1,$ then $u_{\mu_1}$ is a
super-solution of \eq{P} with $\mu=\mu_2.$ It is easy to prove
that $\zeta\leq u_{\mu_1}$ in $\Omega$ and by virtue of Lemma
\ref{sub-sup} we conclude that the problem \eq{P} with $\mu=\mu_2$
has at least one solution.

Thus we have proved $[0,\mu^*)\subseteq A.$ Next we show $\mu^*<+\infty.$\\
Since $\lim_{s\searrow 0}g(s)=+\infty,$ we can choose $s_0>0$ such
that $g(s)>\bar \sigma$ for all $s\leq s_0.$ Let
$$\di \mu_0=\frac{\bar \sigma}
{\min_{x\in\overline\Omega}f(x,s_0)}.$$ Using the monotonicity of
$f$ with respect to the second argument, the above relations
yields
$$\di g(s)+\mu f(x,s)\geq \bar\sigma,\quad \mbox{ for all }\;
(x,s)\in\overline\Omega\times (0,\infty)\;\,\mbox{ and }\;
\mu>\mu_0.$$ If \eq{P} has a solution for $\mu>\mu_0,$ this would
be a super-solution of the problem
\begin{equation}\label {treipatru}
\left\{\begin{tabular}{ll}
$-\Delta u=|\nabla u|^p+\bar\sigma$ \quad & $\mbox{\rm in}\ \Omega,$\\
$u=0$ \quad & $\mbox{\rm on}\ \partial\Omega.$\\
\end{tabular} \right.
\end{equation}
Since 0 is a sub-solution, we deduce that \eq{treipatru} has at
least one solution. According to Lemma \ref{l4}, this is a
contradiction. Hence $\mu^*\leq \mu_0<+\infty.$ This concludes the
proof of Theorem \ref{th1}.
\endproof

\section{Proof of Theorem \ref{th2}}
(i) We fix $p\in(0,1]$ and define
$$q=q(p)=\left\{\begin{tabular}{ll}
$p+1$\quad&${\rm if\; } 0<p<1,$\\
$3/2$\quad&${\rm if\; } p=1.$\\
\end{tabular}\right.$$
Consider the problem
\neweq{treisup}
\left\{\begin{tabular}{ll}
$\di -\Delta u=g(u)+1+|\nabla u|^{q}+\mu f(x,u)$&$\mbox {\rm in}\ \Omega,$\\
$\di u>0$&$\mbox{\rm in}\ \Omega,$\\
$\di u=0$&$\mbox{\rm on}\ \partial\Omega.$\\
\end{tabular}\right.
\endeq
Since $s^p\leq s^q+1,$ for all $s\geq 0,$ we deduce that any
solution of \eq{treisup} is a super-solution of \eq{P}.
Furthermore, taking into account the fact that $1<q<2,$ we can
apply Theorem \ref{th1}(ii) in order to get that \eq{treisup} has
at least one solution if $\mu$ is small enough. Thus, by Lemma
\ref{sub-sup} we deduce that \eq{P} has at least one classical
solution. Following the method used in the proof of Theorem
\ref{th1}, we set
$$\di A=\{\mu\geq 0;\;\mbox{ problem }\eq{P}\mbox{ has at least one solution}\}$$
and let $\mu^*=\sup A.$ With the same arguments we prove that
$[0,\mu^*)\subseteq A.$ It remains only to show that
$\mu^*<+\infty.$
\medskip

Let us assume first that $f$ satisfies $(f1).$ Since
$\lim_{s\searrow 0}g(s)=+\infty,$ we can choose
$\mu_0>\frac{2\la_1}{c}\,$ such that $\frac{1}{2}\mu_0cs+g(s)\geq
1$ for all $s>0.$ Then
$$\di g(s)+\mu f(x,s)\geq \la_1 s+1,\quad\mbox{ for all }\;
(x,s)\in\overline\Omega\times(0,\infty)\;\mbox{ and }\; \mu\geq
\mu_0.$$ By virtue of Lemma \ref{l01} we obtain that \eq{P} has no
classical solutions if $\mu\geq\mu_0,$ so $\mu^*$ is finite.

Assume now that $f$ satisfies $(f2).$ Since $\lim_{s\searrow
0}g(s)=+\infty,$ there exists $s_0>0$ such that
\neweq{f21}
\di g(s)\geq \la_1(s+1)\quad\mbox{ for all }\; 0<s<s_0.
\endeq
On the other hand, the assumption $(f2)$ and the fact that
$\Omega$ is bounded implies that the mapping
$$\di(0,\infty)\ni s\longmapsto
\frac{\min_{x\in\overline\Omega}f(x,s)}{s+1}$$ is nondecreasing,
so we can choose $\tilde\mu>0$ with the property
\neweq{f22}
\tilde\mu\cdot \frac{\min_{x\in\overline\Omega}f(x,s)}{s+1}\geq
\la_1\quad \mbox{ for all }\;s\geq s_0.
\endeq
Now \eq{f21} combined with \eq{f22} yields
$$\di g(s)+\mu f(x,s)\geq \la_1 (s+1),\quad\mbox{ for all }\;(x,s)\in
\overline\Omega\times(0,\infty)\;\mbox { and
}\;\mu\geq\tilde\mu.$$ Using Lemma \ref{l01}, we deduce that
\eq{P} has no solutions if
$\mu>\tilde\mu,$ that is, $\mu^*$ is finite.\\
The first part in Theorem \ref{th2} is therefore established.
\medskip

\noindent (ii) The strategy is to find a super-solution $\overline
u_\mu\in C^2(\Omega)\cap C(\overline\Omega)$ of \eq{P} such that
$\zeta\leq \overline u_\mu$ in $\Omega.$ To this aim, let $h\in
C^2(0,\eta]\cap C[0,\eta]$ be such that
\begin{equation}\label{dunu}
\left\{\begin{tabular}{ll}
$h''(t)=-g(h(t)),\quad \mbox{ for all } 0<t<\eta,$\\
$h(0)=0,$\\
$h>0\quad\mbox{\rm in}\  (0,\eta].$
\end{tabular} \right.
\end{equation}
The existence of $h$ follows by classical arguments of ODE. Since
$h$ is concave, there exists $h'(0+)\in$$(0,+\infty].$ By taking
$\eta>0$ small enough, we can assume that $h'>0$ in $(0,\eta],$ so
$h$ is increasing on $[0,\eta].$
\begin{lem}\label{lh}
{\rm (i)} $h\in C^1[0,\eta]$ if and only if $\di\int_0^1g(s)ds<+\infty;$\\
{\rm (ii)} If $0<p\leq 2,$ then there exist $c_1,c_2>0$ such that
$$\di (h')^p(t)\leq c_1g(h(t))+c_2,\quad\mbox{ for all }\;0<t<\eta.$$
\end{lem}
\noindent\proof (i) Multiplying by $h'$ in \eq{dunu} and then
integrating on $[t,\eta],$ $0<t<\eta,$ we get
\neweq{hprim}
\di (h')^2(t)-(h')^2(\eta)=2\int_t^\eta
g(h(s))h'(s)ds=2\int_{h(t)}^{h(\eta)} g(\tau)d\tau.
\endeq
This gives
\neweq{hpr}
\di (h')^2(t)= 2G(h(t))+(h')^2(\eta)\quad\mbox{ for all }\;
0<t<\eta,
\endeq
where $G(t)=\di\int_t^{h(\eta)}g(s)ds.$ From \eq{hpr} we deduce
that $h'(0+)$ is finite
if and only if $G(0+)$ is finite, so (i) follows.\\
(ii) Let $p\in(0,2].$ Taking into account the fact that $g$ is
nonincreasing, the inequality \eq{hpr} leads to
\neweq{ddoibis}
\di (h')^2(t)\leq 2h(\eta)g(h(t))+(h')^2(\eta),\quad\mbox{ for all
}0<t<\eta.
\endeq
Since $s^p\leq s^2+1,$ for all $s\geq 0,$ from \eq{ddoibis} we
have
\neweq{ddoi}
\di (h')^p(t)\leq c_1g(h(t))+c_2,\quad\mbox{ for all }0<t<\eta
\endeq
where $c_1=2h(\eta)$ and $c_2=(h')^2(\eta)+1.$ This completes the
proof of our Lemma.
\endproof
\medskip

\noindent{ \bf Proof of Theorem \ref{th2} completed.} Let
$p\in(0,1)$ and $\mu\geq 0$ be fixed. We also fix $c>0$ such that
$c\|\varphi_1\|_\infty<\eta.$ By Hopf's maximum principle, there
exist $\delta>0$ small enough and $\theta_1>0$ such that
\neweq{dtrei}
\di |\nabla \varphi_1|>\theta_1 \quad\mbox{ in }\;\Omega_\delta,
\endeq
where $\Omega_\delta:=\{x\in\Omega;\,\mbox{dist}(x,\partial\Omega)
\leq \delta\}.$\\
Moreover, since $\lim_{s\searrow 0}g(h(s))=+\infty,$ we can pick
$\delta$ with the property
\neweq{dpatru}
\di (c\theta_1)^2g(h(c\varphi_1))-3\mu
f(x,h(c\varphi_1))>0\quad\mbox{ in }\;\Omega_\delta.
\endeq
Let $\di
\theta_2:=\inf\limits_{\Omega\setminus\Omega_\delta}\varphi_1>0.$
We choose $M>1$ with
\neweq{dcinci}
M(c\theta_1)^2>3,
\endeq
\neweq{dsase}
Mc\la_1\theta_2h'(c\|\varphi_1\|_{\infty})>3g(h(c\theta_2)).
\endeq
Since $p<1,$ we also may assume
\neweq{dnouap}
(Mc)^{1-p}\la_1(h')^{1-p}(c\|\varphi_1\|_{\infty})\geq
3\|\nabla\varphi_1\|_\infty^p.
\endeq
On the othe hand, by Lemma \ref{lh}(ii) we can choose $M>1$ such
that
\neweq{dsaptebis}
3(h'(c\varphi_1))^p\leq M^{1-p}(c\theta_1)^{2-p}g(h(c\varphi_1))
\quad\mbox{ in }\;\Omega_\delta.
\endeq
The assumption $(f4)$ yields
$$\di
\lim_{s\ri\infty}\frac{3\mu
f(x,sh(c\|\varphi_1\|_\infty))}{sh(c\|\varphi_1\|_\infty)}=0.$$ So
we can choose $M>1$ large enough such that
$$\di \frac{3\mu f(x,Mh(c\|\varphi_1\|_\infty))}{Mh(c\|\varphi_1\|_\infty)}<
\frac{c\la_1\theta_2h'(c\|\varphi_1\|_\infty)}{h(c\|\varphi_1\|_\infty)},$$
uniformly in $\Omega.$ This leads us to
\neweq{dsaptep}
3\mu
f(x,Mh(c\|\varphi_1\|_\infty))<Mc\la_1\theta_2h'(c\|\varphi_1\|_{\infty}),
\quad \mbox{ for all }\;x\in\Omega.
\endeq

For $M$ satisfying \eq{dcinci}-\eq{dsaptep}, we prove that
$\overline u_\mu=Mh(c\varphi_1)$ is a super-solution of \eq{P}. We
have
\neweq{calculp}
\di -\Delta \overline
u_\la=Mc^2g(h(c\varphi_1))|\nabla\varphi_1|^2+
Mc\la_1\varphi_1h'(c\varphi_1) \quad\mbox{ in }\;\Omega.
\endeq
First we prove that
\neweq{dzecep}
Mc^2g(h(c\varphi_1))|\nabla\varphi_1|^2\geq g(\overline
u_\mu)+|\nabla \overline u_\mu|^p+\mu f(x,\overline
u_\mu)\quad\mbox{ in }\;\Omega_\delta.
\endeq
From \eq{dtrei} and \eq{dcinci} we get
\neweq{dunspep}
\di \frac{1}{3}Mc^2g(h(c\varphi_1))|\nabla\varphi_1|^2\geq
g(h(c\varphi_1))\geq g(Mh(c\varphi_1))=g(\overline u_\mu)
\quad\mbox{ in }\;\Omega_\delta.
\endeq
By \eq{dtrei} and \eq{dsaptebis} we also have
\neweq{ddoispep}
\di \frac{1}{3}Mc^2g(h(c\varphi_1))|\nabla\varphi_1|^2\geq
(Mc)^p(h')^p(c\varphi_1))|\nabla\varphi_1|^p=|\nabla \overline
u_\mu|^p \quad\mbox{ in }\;\Omega_\delta.
\endeq
The assumption $(f3)$ and \eq{dpatru} produce
\neweq{dtreispep}
 \frac{1}{3}Mc^2g(h(c\varphi_1))|\nabla\varphi_1|^2\geq
\mu Mf(x,h(c\varphi_1))\geq \mu f(x,Mh(c\varphi_1)) \quad\mbox{ in
}\;\Omega_\delta.
\endeq
Now, by \eq{dunspep}, \eq{ddoispep} and \eq{dtreispep} we conclude
that \eq{dzecep} is fulfilled.
\medskip

Next we prove
\neweq{dpaispep}
Mc\la_1\varphi_1h'(c\varphi_1)\geq g(\overline u_\mu)+|\nabla
\overline u_\mu|^p+\mu f(x,\overline u_\mu)\quad\mbox{ in }\;
\Omega\setminus\Omega_\delta.
\endeq
From \eq{dsase} we obtain
\neweq{dcincispep}
\di \frac{1}{3}Mc\la_1\varphi_1h'(c\varphi_1)\geq
g(h(c\varphi_1))\geq g(Mh(c\varphi_1))=g(\overline
u_\mu)\quad\mbox{ in }\; \Omega\setminus\Omega_\delta.
\endeq
From \eq{dnouap} we get
\neweq{dsaispep}
\di \frac{1}{3}Mc\la_1\varphi_1h'(c\varphi_1)\geq
(Mc)^p(h')^p(c\varphi_1)|\nabla\varphi_1|^p=|\nabla \overline
u_\mu|^p \quad\mbox{ in }\; \Omega\setminus\Omega_\delta.
\endeq
By \eq{dsaptep} we deduce
\neweq{dsaptispep}
\di \frac{1}{3}Mc\la_1\varphi_1h'(c\varphi_1)\geq \mu
f(x,Mh(c\varphi_1))= \mu f(x,\overline u_\mu) \quad\mbox{ in }\;
\Omega\setminus\Omega_\delta.
\endeq
Obviously, \eq{dpaispep} follows now by
\eq{dcincispep}, \eq{dsaispep} and \eq{dsaptispep}.\\
Combining \eq{calculp} with \eq{dzecep} and \eq{dpaispep} we find
that $\overline u_\mu$ is a super-solution of \eq{P}. Moreover,
$\zeta\leq \overline u_\mu$ in $\Omega.$ Applying Lemma
\ref{sub-sup}, we deduce that \eq{P} has at least one solution for
all $\mu\geq 0.$ This finishes the proof of Theorem \ref{th2}.
\endproof

\section{Proof of Theorem \ref{th3}}
The proof case relies on the same arguments used in the proof of
Theorem \ref{th2}. In fact, the main point is to find a
super-solution $\overline u_\la\in
C^2(\Omega)\cap(\overline\Omega)$ of \eq{P}, while $\zeta$ defined
in \eq{treidoi} is a sub-solution. Since $g$ is nonincreasing, the
inequality $\zeta\leq \overline u_\la$ in $\Omega$ can be proved
easily and the existence of solutions to \eq{P} follows by Lemma
\ref{sub-sup}.
\medskip

Define $c,\delta$ and $\theta_1,\theta_2$ as in the proof of
Theorem \ref{th2}. Let M satisfying \eq{dcinci} and \eq{dsase}.
Since $g(h(s))\ri+\infty$ as $s\searrow 0,$ we can choose
$\delta>0$ such that
\neweq{dpa}
\di (c\theta_1)^2g(h(c\varphi_1))-3f(x,h(c\varphi_1))>0\quad\mbox{
in }\;\Omega_\delta.
\endeq

The assumption $(f4)$ produces
$$\di
\lim_{s\ri\infty}\frac{f(x,sh(c\|\varphi_1\|_\infty))}{sh(c\|\varphi_1\|_\infty)}=0,
\quad\mbox{ uniformly for }\;\,x\in\overline\Omega.$$ Thus, we can
take $M>3$ large enough, such that
$$\di \frac{f(x,Mh(c\|\varphi_1\|_\infty))}{Mh(c\|\varphi_1\|_\infty)}<
\frac{c\la_1\theta_2h'(c\|\varphi_1\|_\infty)}{3h(c\|\varphi_1\|_\infty)}.$$
The above relation yields
\neweq{dsapte}
3f(x,Mh(c\|\varphi_1\|_\infty))<
Mc\la_1\theta_2h'(c\|\varphi_1\|_{\infty}),\quad\mbox{for all }
\;\;x\in\overline\Omega.
\endeq
Using Lemma \ref{lh}(ii) we can take $\la>0$ small enough such
that the following inequalities hold
\neweq{dopt}
3\la M^{p-1}(h')^p(c\varphi_1)\leq
g(h(c\varphi_1))(c\theta_1)^{2-p}\quad\mbox{\ in } \;\Omega_\delta
\endeq
\neweq{dnoua}
\la_1\theta_2 h'(c\|\varphi_1\|_\infty)>
 3\la (Mc)^{p-1}
(h')^p(c\theta_2)\|\nabla\varphi_1\|_\infty^p.
\endeq
For $M$ and $\la$ satisfying \eq{dcinci}-\eq{dsase} and
\eq{dpa}-\eq{dnoua}, we claim that $\overline
u_\la=Mh(c\varphi_1)$ is a super-solution of \eq{P}. First we have
\neweq{calcul}
\di -\Delta \overline
u_\la=Mc^2g(h(c\varphi_1))|\nabla\varphi_1|^2+
Mc\la_1\varphi_1h'(c\varphi_1) \quad\mbox{ in }\;\Omega.
\endeq
Arguing as in the proof of Theorem \ref{th2}, from \eq{dtrei},
\eq{dcinci}, \eq{dpa}, \eq{dopt} and the assumption $(f3)$ we
obtain
\neweq{dzece}
Mc^2g(h(c\varphi_1))|\nabla\varphi_1|^2\geq g(\overline
u_\la)+\la|\nabla \overline u_\la|^p+f(x,\overline
u_\la)\quad\mbox{ in }\;\Omega_\delta.
\endeq
On the other hand, \eq{dsase}, \eq{dsapte} and \eq{dnoua} gives
\neweq{dpaispe}
Mc\la_1\varphi_1h'(c\varphi_1)\geq g(\overline u_\la)+\la|\nabla
\overline u_\la|^p+f(x,\overline u_\la)\quad\mbox{ in }\;
\Omega\setminus\Omega_\delta.
\endeq
Using now \eq{calcul} and \eq{dzece}-\eq{dpaispe} we find that
$\overline u_\la$
is a super-solution of \eq{P} so our claim follows.\\
As we have already argued at the beginning of this case, we easily
get that $\zeta\leq\overline u_\la$ in $\Omega$ and by Lemma
\ref{sub-sup} we deduce that problem \eq{P} has at least one
solution
if $\la>0$ is sufficiently small.\\
Set
$$ \di A=\{\;\la\geq 0; \mbox{ problem }\eq{P}\mbox{ has at least one classical solution}\}.$$
From the above arguments, $A$ is nonempty. Let $\la^*=\sup A.$
First we claim that if $\la\in A,$ then $[0,\la)\subseteq A.$ For
this purpose, let $\la_1\in A$ and $0\leq \la_2<\la_1.$ If
$u_{\la_1}$ is a solution of \eq{P} with $\la=\la_1,$ then
$u_{\la_1}$ is a super-solution for \eq{P} with $\la=\la_2$ while
$\zeta$ defined in \eq{treidoi} is a sub-solution. Using Lemma
\ref{sub-sup} once more, we have that \eq{P} with $\la=\la_2$ has
at least one classical solution. This proves the claim. Since
$\la\in A$ was arbitrary chosen, we conclude that
$[0,\la^*)\subset A.$
\medskip

Let us assume now $p\in(1,2].$ We prove that $\la^*<+\infty.$ Set
$$\di m:=\inf_{(x,s)\in\overline\Omega\times(0,\infty)}\Big(g(s)+f(x,s)\Big).$$
Since $\lim_{s\searrow 0}g(s)=+\infty$ and the mapping $\di
(0,\infty)\ni s\longmapsto \min_{x\in\overline\Omega}f(x,s)$ is
positive and nondecreasing, we deduce that $m$ is a positive real
number. Let $\la>0$ be such that \eq{P} has a solution $u_\la.$ If
$v=\la^{1/(p-1)}u_\la,$ then $v$ verifies
\neweq{doptispe}
\left\{\begin{tabular}{ll} $\di -\Delta v\geq |\nabla
v|^p+\la^{1/(p-1)}m$ \quad & $\mbox{\rm in }\Omega,$\\
$v>0$ \quad & $\mbox{\rm in }\Omega,$\\
$v=0$ \quad & $\mbox{\rm on }\, \partial\Omega.$\\
\end{tabular} \right.
\endeq
It follows that $v$ is a super-solution of \eq{doitrei} for $\di
\sigma=\la^{1/(p-1)}m.$ Since 0 is a sub-solution, we obtain that
\eq{doitrei} has at least one classical solution for $\sigma$
defined above. According to Lemma \ref{l4}, we have $\sigma\leq
\bar\sigma,$ and so $\di
\la\leq\left(\frac{\bar\sigma}{m}\right)^{p-1}.$ This means that
$\la^*$ is finite.
\medskip

Assume now $p\in(0,1)$ and let us prove that $\la^*=+\infty.$
Recall that $\zeta$ defined in \eq{treidoi} is a sub-solution. To
get a super-solution, we proceed in the same manner. Fix $\la>0.$
Since $p<1$ we can find $M>1$ large enough such that
\eq{dcinci}-\eq{dsase} and \eq{dsapte}-\eq{dnoua} hold.
From now on, we follow the same steps as above.\\
The proof of Theorem \ref{th3} is now complete.
\endproof
\noindent {\bf Remark.} If $\di\int _0^1g(s)ds<\infty,$ then the
above method can be applied in order to extend the study of \eq{P}
to the case $\mu=1$ and $p>2.$ Indeed, by Lemma \ref{lh}(i) it
follows $h\in C^1[0,\eta].$ Using this fact, we can choose
$c_1,c_2>0$ large enough such that the conclusion of Lemma
\ref{lh}(ii) holds. Repeating the above arguments we prove that if
$p>2$ then there exists a real number $\la^*>0$ such that \eq{P}
has at least one solution if $\la<\la^*$ and no solutions exist if
$\la>\la^*.$

\section{Proof of Theorem \ref{th4}}
(i) If $\la=0,$ the existence of the solution follows by using
Lemma \ref{l2}.
 Next we assume that $\la>0$
and let us fix $\mu\geq 0.$ With the change of variable $v=e^{\la
u}-1,$ the problem \eq{P} becomes
\neweq{Q}
\left\{\begin{tabular}{ll} $-\Delta v=\Phi_\la(v)$ \quad &
$\mbox{\rm
in } \Omega,$\\
$v>0$ \quad & $\mbox{\rm in } \Omega,$\\
$v=0$ \quad & $\mbox{\rm on } \partial\Omega,$\\
\end{tabular} \right.
\endeq
where
$$\di
\Phi_\la(s)=\la(s+1)g\left(\frac{1}{\la}\ln(s+1)\right)+\la\mu
(s+1),$$ for all $s\in (0,\infty).$ Obviously $\Phi_\la$ is not
monotone but we still have that the mapping $\di (0,\infty)\ni
s\mapsto \frac{\Phi_\la(s)}{s}\,$ is decreasing for all $\la>0$
and
$$\di \lim _{s\rightarrow+\infty}
\frac{\Phi_\la(s)}{s}=\la(a+\mu) \quad\mbox{ and }\quad
\lim_{s\searrow
0}\frac{\Phi_\la(s)}{s}=+\infty,$$ uniformly for $\la>0.$\\
We first remark that $\Phi_{\la}$ satisfies the hypotheses in
Lemma \ref{l2} provided $\la(a+\mu)<\la_1.$ Hence \eq{Q} has at
least one solution.

\noindent On the other hand, since $g\geq a$ on $(0,\infty),$ we
get
\neweq{Phi}
\di \Phi_\la(s)\geq\la(a+\mu)(s+1),\quad\mbox{ for all }\;
\la,s\in(0,\infty).
\endeq
Using now Lemma \ref{l01} we deduce that \eq{Q} has no solutions
if $\la(a+\mu)\geq \la_1.$ The proof of the first part in Theorem
\ref{th4} is therefore complete.
\medskip

\noindent(ii) We split the proof into several steps.

\noindent {\sc Step 1.} {\bf Existence of solutions.}\\
This follows directly from (i).
\medskip

\noindent {\sc Step 2.} {\bf Uniqueness of the solution.}\\
Fix $\la\geq 0.$ Let $u_1$ and $u_2$ be two classical solutions of
\eq{P} with $\la<\la^*.$ We show that $u_1\leq u_2$ in $\Omega.$
Supposing the contrary, we deduce that
$\max\limits_{\overline\Omega}\{u_1-u_2\}>0$ is achieved in a
point $x_0\in \Omega.$ This yields $\nabla (u_1-u_2)(x_0)=0$ and
$$\di 0\leq -\Delta
(u_1-u_2)(x_0)=g(u_1(x_0))-g(u_2(x_0))<0,$$ a contradiction. We
conclude that $u_1\leq u_2$ in $\Omega;$ similarly $u_2\leq u_1.$
Therefore $u_1= u_2$ in $\Omega$ and the uniqueness is proved.
\medskip

\noindent {\sc Step 3.} {\bf Dependence on $\la$.}\\
Fix $0\leq \la_1<\la_2<\la^*$ and let $u_{\la_1},$ $u_{\la_2}$ be
the unique solutions of \eq{P} with $\la=\la_1$ and $\la=\la_2$
respectively. If $\{x\in\Omega;u_{\la_1}>u_{\la_2}\}$ is nonempty,
then $\max\limits_{\overline\Omega}\{u_{\la_1}-u_{\la_2}\}>0$ is
achieved in $\Omega.$ At that point, say $\bar x,$ we have $\nabla
(u_{\la_1}-u_{\la_2})({\bar x})=0$ and
$$0\leq -\Delta(u_{\la_1}-u_{\la_2})(\bar x)=g(u_{\la_1}(\bar x))-g(u_{\la_2}(\bar x))
+(\la_1-\la_2)|\nabla u_{\la_1}|^p(\bar x)<0,$$
which is a contradiction.\\
Hence $u_{\la_1}\leq u_{\la_2}$ in $\overline\Omega.$ The maximum
principle also gives $u_{\la_1}< u_{\la_2}$ in $\Omega.$
\medskip

\noindent {\sc Step 4.} {\bf Regularity.}\\
We fix $0<\la<\la^*,$ $\mu>0$ and assume that $\limsup_{s\searrow
0}s^{\alpha}g(s)<+\infty.$ This means that $g(s)\leq cs^{-\alpha}$
in a small positive neighborhood of the origin. To prove the
regularity, we will use again the change of variable $v=e^{\la
u}-1.$ Thus, if $u_\la$ is the unique solution of \eq{P}, then
$v_\la=e^{\la u_\la}-1$ is the unique solution of \eq{Q}. Since
$\di \lim_{s\searrow 0}\frac{e^{\la s}-1}{s}=\la,$ we conclude
that (ii1) and (ii2) in Theorem \ref{th4} are established if we
prove

${\rm (a)}\quad \tilde c_1\,{\rm dist}(x,\partial \Omega) \leq
v_\la(x)\leq \tilde c_2\,{\rm dist} (x,\partial \Omega)$ in
$\Omega,$ for some positive constants $\tilde c_1,\tilde c_2>0.$

${\rm (b)}\quad v_\la\in C^{1,1-\alpha}(\overline \Omega).$\\
{\it Proof of} (a). By the monotonicity of $g$ and the fact that
$g(s)\leq cs^{-\alpha}$ near the origin, we deduce the existence
of $A,B,C>0$ such that
\neweq{ctrei}
\di \Phi_\la(s)\leq As+Bs^{-\alpha}+C,\quad \mbox{ for
all}\;0<\la<\la^*\mbox { and }\,s>0.
\endeq
Let us fix $m>0$ such that $m\la_1\|\varphi_1\|_{\infty}<\la\mu.$
Combining this with \eq{Phi} we deduce
\neweq{sssu}
\di -\Delta(v_\la-m\varphi_1)=\Phi_\la(v_\la)-m\la_1\varphi_1\geq
\la\mu-m\la_1\varphi_1\geq 0
\endeq
in $\Omega.$ Since $v_\la-m\varphi_1=0$ on $\partial\Omega,$ we
conclude
\neweq{fi2}
v_\la\geq m\varphi_1\quad\mbox{ in }\;\Omega.
\endeq
Now, \eq{fi2} and \eq{fi} imply $v_\la\geq \tilde c_1\,{\rm dist}
(x,\partial \Omega)$ in $\Omega,$ for some positive constant
$\tilde c_1>0.$ The first inequality in the statement of (a) is
therefore established. For the second one, we apply an idea found
in Gui and Lin \cite{gl}. Using \eq{fi2} and the estimate
\eq{ctrei}, by virtue of Lemma \ref{l1} we deduce
$\Phi_\la(v_\la)\in L^1(\Omega),$ that is, $\Delta v_\la\in
L^1(\Omega).$

Using the smoothness of $\partial\Omega,$ we can find
$\delta\in(0,1)$ such that for all
$x_0\in\Omega_{\delta}:=\{x\in\Omega\,;\,{\rm
dist}(x,\partial\Omega)\leq \delta\},$ there exists
$y\in\RR^N\setminus\overline\Omega$ with ${\rm
dist}(y,\partial\Omega)=\delta$ and ${\rm
dist}(x_0,\partial\Omega)=|x_0-y|-\delta.$

Let $K>1$ be such that diam$(\Omega)<(K-1)\delta$ and let $\xi$ be
the unique solution of the Dirichlet problem
\neweq{bdoi}
\left\{\begin{tabular}{ll} $-\Delta \xi=\Phi_\la(\xi)$ \quad &
${\rm
in}\ B_K(0)\setminus B_1(0),$\\
$\xi>0$ \quad & ${\rm in}\ B_K(0)\setminus B_1(0),$\\
$\xi=0$ \quad & ${\rm on}\ \partial(B_K(0)\setminus B_1(0)).$\\
\end{tabular} \right.
\end{equation}
where $B_r(0)$ denotes the open ball in $\RR^N$ of radius $r$ and
centered at the origin. By uniqueness, $\xi$ is radially
symmetric. Hence $\xi(x)=\tilde \xi(|x|)$ and
\neweq{btrei}
\left\{\begin{tabular}{ll} $\di \tilde \xi''+\frac{N-1}{r}\tilde
\xi'+\Phi_\la(\tilde\xi)=0$
\quad & ${\rm in}\ (1,K),$\\
$\tilde \xi>0$ \quad & ${\rm in}\ (1,K),$\\
$\tilde \xi(1)=\tilde \xi(K)=0.$ \quad & \\
\end{tabular} \right.
\end{equation}
Integrating in \eq{btrei} we have
$$\begin{aligned}
\di \tilde \xi'(t)& =\tilde \xi'(a)a^{N-1}t^{1-N}-t^{1-N}\int_a^t
r^{N-1} \Phi_\la(\tilde \xi(r))dr\\
& =\tilde \xi'(b)b^{N-1}t^{1-N}+t^{1-N}\int_t^b
r^{N-1} \Phi_\la(\tilde \xi(r))dr,\\
\end{aligned}$$
where $1<a<t<b<K.$ With the same arguments as above we have
$\Phi_\la(\tilde \xi)\in L^1(1,K)$ which implies that both $\tilde
\xi(1)$ and $\tilde \xi(K)$ are finite. Hence $\tilde \xi\in
C^2(1,K)\cap C^1[1,K].$ Furthermore,
\neweq{bpatru}
\xi(x)\leq \tilde C\min\{K-|x|,|x|-1\}, \quad\mbox{ for any }
\;\;x\in B_K(0)\setminus B_1(0).
\endeq
Let us fix $x_0\in\Omega_{\delta}.$ Then we can find
$y_0\in\RR^N\setminus\overline\Omega$ with ${\rm
dist}(y_0,\partial\Omega)=\delta$ and ${\rm
dist}(x_0,\partial\Omega)=|x_0-y|-\delta.$ Thus, $\Omega\subset
B_{K\delta}(y_0)\setminus B_{\delta}(y_0).$ Define $\di \overline
v(x)=\xi\left(\frac{x-y_0}{\delta}\right),$ for all
$x\in\overline\Omega.$ We show that $\overline v$ is a
super-solution of \eq{Q}. Indeed, for all $x\in\Omega$ we have
$$\begin{tabular}{rl}
$\di\Delta \overline v+\Phi_\la(\overline v)$ &
$\di=\frac{1}{\delta^2}\left(\tilde \xi''+
\frac{N-1}{r}\tilde \xi'\right)+\Phi_\la(\tilde \xi)$\\
&$\di \leq\frac{1}{\delta^2}\left(\tilde \xi''+
\frac{N-1}{r}\tilde \xi'+\Phi_\la(\tilde \xi)\right)$\\
&$=0,$\\
\end{tabular}$$
where $\di r=\frac{|x-y_0|}{\delta}.$ We have obtained that
$$\di \Delta \overline v+\Phi_\la(\overline v)\leq 0\leq \Delta v_\la+\Phi_\la(v_\la)
\quad\mbox{ in }\;\Omega,$$
$$\di \overline v,v_\la>0\;\;\mbox{ in }\;\Omega,\;\overline v=v_\la
\;\;\mbox{ on }\;\partial\Omega $$
$$\Delta v_\la\in L^1(\Omega).$$
By Lemma \ref{l3} we get $\di v_{\la}\leq \overline v$ in
$\Omega.$ Combining this with \eq{bpatru} we obtain
$$\di v_{\la}(x_0)\leq \overline v(x_0)\leq
\tilde C\min\{K-\frac{|x_0-y_0|}{\delta},
\frac{|x_0-y_0|}{\delta}-1\}\leq\frac{\tilde C}{\delta}{\rm
dist}(x_0,\partial\Omega).$$ Hence $v_{\la}\leq \frac{\tilde
C}{\delta}{\rm dist}(x,\partial\Omega)$ in $\Omega_{\delta}$ and
the second inequality in the statement of (a) follows.
\medskip

\noindent {\it Proof of} (b). Let $G$ be the Green's function
associated with the Laplace operator in $\Omega.$ Then, for all
$x\in\Omega$ we have
$$\di v_{\la}(x)=-\int_{\Omega} G(x,y)\Phi_\la(v_{\la}(y))dy$$
and
$$\di \nabla v_{\la}(x)=-\int_{\Omega} G_x(x,y)\Phi_\la(v_{\la}(y))dy.$$
If $x_1,x_2\in\Omega,$ using \eq{ctrei} we obtain
$$\begin{tabular}{ll}
$\di |\nabla v_{\la}(x_1)-\nabla v_{\la}(x_2)|$&$\di
\leq \int_{\Omega} |G_x(x_1,y)-G_x(x_2,y)|\cdot(Av_\la+C)dy$\\
&$\di \;\;\,+B\int_{\Omega}
|G_x(x_1,y)-G_x(x_2,y)|\cdot v_{\la}^{-\alpha}(y)dy.$\\
\end{tabular}$$
Now, taking into account that $v_{\la}\in C(\overline\Omega),$ by
the standard regularity theory (see \cite{gt}) we get
$$\di \int_{\Omega} |G_x(x_1,y)-G_x(x_2,y)|\cdot
(Av_\la+C)dy\leq \tilde c_1|x_1-x_2|.$$ On the other hand, with
the same proof as in \cite[Theorem 1]{gl}, we deduce
$$\di \int_{\Omega} |G_x(x_1,y)-G_x(x_2,y)|\cdot
v_{\la}^{-\alpha}(y)\leq \tilde c_2 |x_1-x_2|^{1-\alpha}.$$ The
above inequalities imply $u_{\la}\in C^2(\Omega)\cap
C^{1,1-\alpha}({\overline{\Omega}}).$

\noindent {\sc Step 5.} {\bf Asymptotic behaviour of the solution.}\\
In order to conclude the asymptotic behaviour for $u_\la,$ it is
enough to show that $\lim_{\la\nearrow \la^*}v_\la=+\infty$ on
compact subsets of $\Omega.$ To this aim, we use some techniques
developed in \cite{mr}. Due to the special character of our
problem, we will be able to show in what follows that, in certain
cases, $L^2$--boundedness implies $H^1_0$--boundedness!

We argue by contradiction. Since $(v_{\la})_{\la<\la^*}$ is a
sequence of nonnegative super-harmonic functions in $\Omega$ then,
by \cite[Theorem 4.1.9]{h}, we can find a subsequence of
$(v_{\la})_{\la<\la^*}$ (still denoted by $(v_{\la})_{\la<\la^*}$
) which converges in $L^1_{\rm loc}(\Omega)$ to some $v^*.$ The
monotony of $v_{\la}$ yields (up to a subsequence)
$v_{\la}\nearrow v^*$ a.e. in $\Omega.$
\medskip

We first show that $(v_{\la})_{\la<\la^*}$ is bounded in
$L^2(\Omega).$ Suppose the contrary. Passing eventually at a
subsequence, we have $v_{\la}=M(\la)w_{\la},$ where
\neweq{cdoi}
M(\la)=||v_{\la}||_{L^2(\Omega)}\ri\infty \;\;\mbox { as
}\;\;\la\nearrow\la^* \quad\mbox{ and }\;\; w_{\la}\in
L^2(\Omega),\;\; ||w_{\la}||_{L^2(\Omega)}=1.
\endeq
Then \eq{ctrei} yields
$$\di \frac{1}{M(\la)}\Phi_\la(v_\la)\ri 0 \quad
\mbox{ in } L^1_{\rm loc}(\Omega) \;\;\mbox { as
}\;\;\la\nearrow\la^*$$ that is,
\neweq{cpatru}
\di -\Delta w_{\la}\ri 0 \quad\mbox{\ii n } L^1_{\rm loc}(\Omega)
\;\;\mbox { as }\;\;\la\nearrow\la^*.\end{equation} By Green's
first identity, we have
\neweq{wwbc}
\int_\Omega \nabla w_\la\cdot\nabla \phi dx=-\int_\Omega\phi\Delta
w_\la dx= -\int_{{\rm Supp }\,\phi}\phi\Delta w_\la dx,\quad\mbox{
for all }\;\phi\in C^{\infty}_0(\Omega).
\endeq
Using \eq{cpatru} we obtain
\neweq{wwbbc}
\begin{tabular}{ll}
$\di \left|\int_{{\rm Supp }\,\phi}\phi\Delta w_\la dx\right|$&
$\di\leq \int_{{\rm Supp }\,\phi}|\phi||\Delta w_\la |dx$\\
&$\di\leq \|\phi\|_{\infty}\int_{{\rm Supp }\,\phi}|\Delta
w_\la|dx\ri 0 \mbox{ as }\;\la\nearrow \la^*.$
\end{tabular}\endeq
Now, \eq{wwbc} and \eq{wwbbc} yield
\neweq{wwbbbc}
\di\int_{\Omega}\nabla w_\la\cdot\nabla\phi dx\ri 0\quad\mbox{ as
}\;\la\nearrow\la^*,\; \mbox{ for all }\;\phi\in
C^\infty_0(\Omega).
\endeq

Recall that $(w_{\la})_{\la<\la^*}$ is bounded in $L^2(\Omega).$
We claim that $(w_{\la})_{\la<\la^*}$ is bounded in
$H^1_0(\Omega).$ Indeed, using \eq{ctrei} and H\"{o}lder's
inequality, we have
$$\begin{tabular}{ll}
$\di \int_{\Omega}|\nabla w_{\la}|^2$&$\di
=-\int_{\Omega}w_{\la}\Delta w_{\la}$\\
&$\di =-\frac{1}{M(\la)}\int_{\Omega} w_{\la}
\Delta u_{\la}=\frac{1}{M(\la)}\int_\Omega w_\la\Phi_\la(v_\la)$\\
&$\di \leq \frac{A}{M(\la)}\int_{\Omega}
w_{\la}v_{\la}+\frac{B}{M(\la)}
\int_{\Omega}w_{\la}v^{-\alpha}_{\la}+\frac{C}{M(\la)}\int_\Omega w_\la$\\
&$\di=A\int_{\Omega}w^2_{\la}+
\frac{B}{M(\la)^{1+\alpha}}\int_{\Omega}w^{1-\alpha}_{\la}+\frac{C}{M(\la)}
\int_{\Omega}w_{\la}$\\
&$\di \leq A+\frac{B}{M(\la)^{1+\alpha}}|\Omega|^{(1+\alpha)/2}+
\frac{C}{M(\la)}|\Omega|^{1/2}.$\\
\end{tabular}$$
From the above estimates, we can easily conclude that
$(w_{\la})_{\la<\la^*}$ is bounded in $H^1_0(\Omega).$ Thus, there
exists $w\in H^1_0(\Omega)$ such that
\neweq{ccinci}
w_{\la}\weak w \quad\mbox{ weakly in
}\;\;H^1_0(\Omega)\end{equation} and
\neweq{csase}
w_{\la}\ri w \quad\mbox{ strongly in
}\;\;L^2(\Omega).\end{equation} Combining \eq{cdoi} and
\eq{csase}, we get $\|w\|_{L^2(\Omega)}=1.$ On the other hand,
from \eq{wwbbbc} and \eq{ccinci} we obtain
$$\di\int_{\Omega}\nabla w\cdot\nabla \phi dx=0,
\quad\mbox{ for all }\;\phi\in C^\infty_0(\Omega).$$ Since $w\in
H^1_0(\Omega),$ using the above relation and the definition of
$H^1_0(\Omega),$ we get $w=0,$ which contradicts the fact that
$\|w\|_{L^2(\Omega)}=1.$ Hence $(v_{\la})_{\la<\la^*}$ is bounded
in $L^2(\Omega).$ As before for $w_\la,$ we can obtain that
$(v_\la)_{\la<\la^*}$ is bounded in $H_0^1(\Omega).$ Then, up to a
subsequence we have
\neweq{wbbbcc}
\begin{tabular}{ll}
$\di v_\la\weak v^*$&$\mbox{ weakly in }\;H^1_0(\Omega)\;\,\mbox{
as }\;\,
\la\nearrow\la^*,$\\
$\di v_\la\ri v^*$&$\mbox{ strongly in }\; L^2(\Omega)\;\,
\mbox{ as }\la\nearrow\la^*,$\\
$\di v_\la\ri v^*$&$\mbox{ a.e. in }\;\Omega\;\,\mbox{ as
}\;\la\nearrow \la^*.$
\end{tabular}\endeq
Now we can proceed to get a contradiction. Multiplying by
$\varphi_1$ in \eq{Q} and then integrating over $\Omega$ we have
\neweq{csapte}
\di -\int_{\Omega}\Delta v_{\la}\varphi_1
dx=\int_{\Omega}\Phi_\la(v_\la)\varphi_1dx \quad\mbox{ for
all}\;0<\la<\la^*.
\end{equation}
Using \eq{Phi} we get
\neweq{copt}
\di \la_1\int_{\Omega}v_{\la}\varphi_1\geq \la(a+\mu)
\int_{\Omega}(v_{\la}+1)\varphi_1dx,\quad\mbox{ for all
}\;0<\la<\la^*.
\endeq
By \eq{wbbbcc} we can use Lebesgue's dominated convergence theorem
in order to pass to the limit with $\la\nearrow\la^*$ in
\eq{copt}. We obtain
\begin{equation}
\di \la_1\int_{\Omega}v^*\varphi_1dx\geq
\la_1\int_{\Omega}(v^*+1)\varphi_1dx,
\end{equation}
which is a contradiction since $\varphi_1>0$ in $\Omega.$ This
contradiction shows that $\di \lim_{\la\nearrow
\la^*}v_{\la}=+\infty,$ uniformly on compact subsets of $\Omega$
which implies $\di \lim_{\la\nearrow \la^*}u_{\la}=+\infty,$
uniformly on compact subsets of $\Omega.$ The proof of Theorem
\ref{th4} is now complete.
\endproof

\end{document}